\documentclass[11pt]{amsart}
\usepackage{amsmath,amsthm}
\usepackage{graphicx}
\usepackage{animate}
\usepackage{subfigure}
\usepackage{mathptmx}
\usepackage{latexsym}
\usepackage[T1]{fontenc}
\usepackage[latin1]{inputenc} 
\usepackage{amssymb}
\usepackage[center]{caption}
\usepackage{hyperref}

\begin{document}

\title[Eta Function Series Partial Sums]{Partial Sums of the Series for the Dirichlet Eta Function, their Peculiar Convergence, the Simple Zeros Conjecture, and the RH.} 

\author{Luca Ghislanzoni}

\dedicatory{luca.ghislanzoni@c-sigma.it}  

\date{July 16, 2025}

\begin{abstract} 
For any $s \in \mathbb{C}$\, with $\Re(s)>0$, denote by $\eta_{n-1}(s)$\, the $(n-1)^{th}$\, partial sum of the Dirichlet series for the eta function $\eta(s)=1-2^{-s}+3^{-s}-\cdots \;$\,, and by $R_n(s)$ the corresponding remainder. Denoting by $u_n(s)$ the segment starting at $\eta_{n-1}(s)$ and ending at $\eta_n(s)$, we first show how, for sufficiently large $n$ values, the circle of diameter $u_{n+2}(s)$ lies strictly inside the circle of diameter $u_n(s)$, to then derive the asymptotic relationship $R_n(s) \sim (-1)^{n-1}/n^s$, as $n \rightarrow \infty$. Denoting by $D=\left\{s \in \mathbb{C}: \; 0< \Re(s) < \frac{1}{2}\right\}$\, the open left half of the critical strip, define for all $s\in D$ the ratio $\chi_n^{\pm}(s)  =  \eta_n(1-s) /  \eta_n(s)$. We then prove that the limit $L(s)=\lim_{N(s)<n\to\infty} \chi_n^{\pm}(s)$\, exists at every point $s$\, of the domain $D$\,. The function $L(s)$\, is continuous on $D$\, if and only if the Riemann Hypothesis is true. Finally, we remark how the asymptotic behavior of $R_n(s)$ can also provide insights substantiating the so called Simple Zeros Conjecture.      
\end{abstract}
\maketitle
\section{Introduction}
\label{sec:1}

On the 11th of August 1859, Bernhard Riemann was appointed member of the Berlin Academy. Respectful of such great honor, Riemann submitted to the Academy his seminal work on the distribution of prime numbers less than a given quantity (Über die Anzahl der Primzahlen unter einer gegebenen Grösse, \cite{01} ). In that paper Riemann formulated a daring hypothesis \cite{02}:
\medskip
\begin{center}
	\textbf{All non-trivial zeros of the zeta function lie on the critical line}
\end{center}

\indent Denoted by $s=\sigma+it$\, a complex number, said Riemann Zeta function, $\zeta(s)$\,, is the function defined as the analytic continuation to all the complex values $s\neq 1$\, of the infinite series:
\begin{equation}
	\zeta(s) = \sum_{n=1}^\infty\frac{1}{n^s} =
	1+\frac{1}{2^s}+\frac{1}{3^s}+\frac{1}{4^s}+\ldots \; ,
	\label{eq_1}  
\end{equation}   
which converges only in the half-plane $ \Re(s)>1 $\,. $\zeta(s)$\, has zeros at $s=-2, \; -4, \; -6, \; \ldots $\,, which are called \textit{trivial zeros} because their existence is easy to prove. The Riemann zeta function features also other zeros, called \textit{non-trivial zeros}, known to lie in the open strip $ {s \in \mathbb{C} : 0<\Re(s)<1} $\,, which is called the \textit{critical strip}. The \textit{critical line} is then defined as the line  ${s \in \mathbb{C} : \Re(s)=\frac{1}{2}} $\,. The Riemann Hypothesis asserts that all \textit{non-trivial zeros} lie on the critical line.\\
\indent Since Riemann's milestone paper the properties of $\zeta(s)$\, have been studied in much depth. Besides Riemann's own results on the distribution of prime numbers, the Riemann Hypothesis has important implications also for a number of problems in physics \cite{03}. From the great wealth of existing literature only few key definitions and results are recalled in this introduction. This work will then concentrate its attention on the behavior of the zeros of $\zeta(s)$\, in the interior of the critical strip.\\
\indent A key role in the analytical continuation to $ \Re(s)>0 $\, of the series (\ref{eq_1}) is played by the Dirichlet Eta function \cite{04}:
\begin{equation}
	\eta(s) = \sum_{n=1}^\infty\frac{(-1)^{n-1}}{n^s} =
	1-\frac{1}{2^s}+\frac{1}{3^s}-\frac{1}{4^s}+-\ldots
	\label{eq_2}
\end{equation}
The above series represents the simplest alternating signs case among Ordinary Dirichlet Series, and it is converging for all $s$\, with $ \Re(s)>0 $\,. Its sum, $\eta(s)$\,, is an analytic function in the corresponding half-plane. Analytic continuation to $ \Re(s)>0 $\, of the infinite series (\ref{eq_1}) can then be obtained by observing that:
\begin{equation}
	\zeta(s) = \frac{\eta(s)}{1-\frac{2}{2^s}} \; ,
	\label{eq_3}
\end{equation}  
which is valid for $ s \neq 1\!+\!n\!\cdot\!\frac{2\pi}{\ln{2}}i\quad(n\in\mathbb{Z}) $\,, values at which the denominator vanishes.\smallskip \\
The right hand side of (\ref{eq_3}) is analytic in the region of interest for this work: the interior of the critical strip. Inside such region, the zeros of the Riemann $ \zeta $\, coincide with the zeros of the $ \eta $\, function. As the infinite sum (\ref{eq_2}) converges readily, it makes it easy to graphically represent the path described by the partial sums.\\
\indent Let us denote the partial sums as $ \eta_{k-1}(s) $\, and the corresponding remainder as $ R_k(s) $\, so defined:
\begin{equation}
	\sum_{n=1}^\infty\frac{(-1)^{n-1}}{n^s} =
	\sum_{n=1}^{k-1}\frac{(-1)^{n-1}}{n^s}+\sum_{n=k}^\infty\frac{(-1)^{n-1}}{n^s}=
	\eta_{k-1}(s) + R_k(s) \; ,
	\label{eq_4}
\end{equation}
\smallskip \indent The terms of the infinite sum (\ref{eq_2}) are complex numbers, which can be represented by vectors of the form $u_n(s)=(-1)^{n-1} \frac{1}{n^{\sigma}}e^{-i t \ln n}$\,. \smallskip The line segments making up the paths graphed in Fig. 1 represent said vectors added up \textit{"tip to tail"}. As it will be better detailed in the proof of \textbf{Theorem 1}, while approaching the point of convergence, $\eta(s)$\,, the path described by the partial sums eventually ends up following a very simply structured \textit{star-shaped} path (depicted in the square inset at the top left of Fig. 1), characterized by angles between consecutive segments being $<\frac{\pi}{2}$\,, getting smaller and smaller as $n$\, grows larger and larger.\smallskip \\
\indent A remarkable functional equation satisfied by $\zeta(s)$\, was originally proposed by Euler in 1749, and later proved by Riemann in his 1859 paper \cite{01} \cite{05} 
\begin{equation}
	\zeta(1-s) = 2 (2\pi)^{-s} \cos\left(\frac{\pi s}{2}\right)\ \Gamma(s)\ \zeta(s) \;\;\;\;\;\;\;\; OR \;\;\;\;\;\;\; \zeta(1-s) = \chi(s) \; \zeta(s)  .
	\label{eq_5}
\end{equation}
\indent When studying the behavior of $\zeta(s)$\, inside the critical strip, a useful implication of (\ref{eq_5}) is that to $\zeta(s)=0$\, it must correspond $\zeta(1-s)=0$\,. Therefore, if $s=\frac{1}{2}-\alpha+ it$\, ($0\leq\alpha<\frac{1}{2}$\,) is a zero of $\zeta(s)$\,, then so it must be for both $1-s=\frac{1}{2}+\alpha - it$\, and its complex conjugate $1-\overline{s}=\frac{1}{2}+ \alpha + it$\,. Thus, non-trivial zeros always occur in groups of two pairs, one pair being the complex conjugate of the other, symmetrically located about the critical line. As this work is concerned with the study of the behavior of the zeros of $\zeta(s)$\  in the interior of the critical strip, and inside such region a useful implication of (\ref{eq_3}) is that said zeros coincide with the zeros of $\eta(s)$\,, one could as well concentrate on the study of the behavior of pairs of critical line symmetrical zeros of the function $\eta(s)$\,, which also satisfies a functional equation, readily obtainable from (\ref{eq_3}) when restricting to the the open critical strip
\begin{equation}
	\eta(1-s) = \frac{1-2^{s}}{1-2^{1-s}} \;2 (2\pi)^{-s} \cos\left(\frac{\pi s}{2}\right)\ \Gamma(s)\ \eta(s) \;\;\;\;\;\;\;\; OR \;\;\;\;\;\;\; \eta(1-s) = \chi^{\pm}(s) \; \eta(s)  .
	\label{eq_6}
\end{equation}
whereby the $^{\pm}$ label is added because this is the functional equation for the alternating sign $\zeta$ function.\\
\indent 
The readily converging sum (\ref{eq_2}) allows to easily visualize, by drawing the corresponding graphs, the paths described by the partial sums of critical line symmetrical arguments. \\
The pattern of convergence of the partial sums, $\eta_N(s)$\,, of the alternating Dirichlet series displays a beautiful and features rich geometrical structure, of which Fig.1 provides only an example (restricted to the first 1000 terms of the sum for $\eta\left(\frac{1}{2} \pm 0.15 + i\,2103.19 \right)$\,). The red path is obtained by drawing the line segments, $u_{N}(s)=(-1)^{N-1}/N^s$ joining partial sum $\eta_{N-1}(s)$\ with  $\eta_{N}(s)$, while the green one is obtained by drawing the corresponding $u_N(1-\bar s)$. The green path represents a sort of image of the red one, in the sense that each green $u_N(1-\bar s)$ segment is by definition parallel to the corresponding red $u_N( s)$ segment, while being scaled down by $N^{2\sigma-1}$. 
\begin{figure}[t]
	\centering
	\includegraphics[width=1\textwidth]{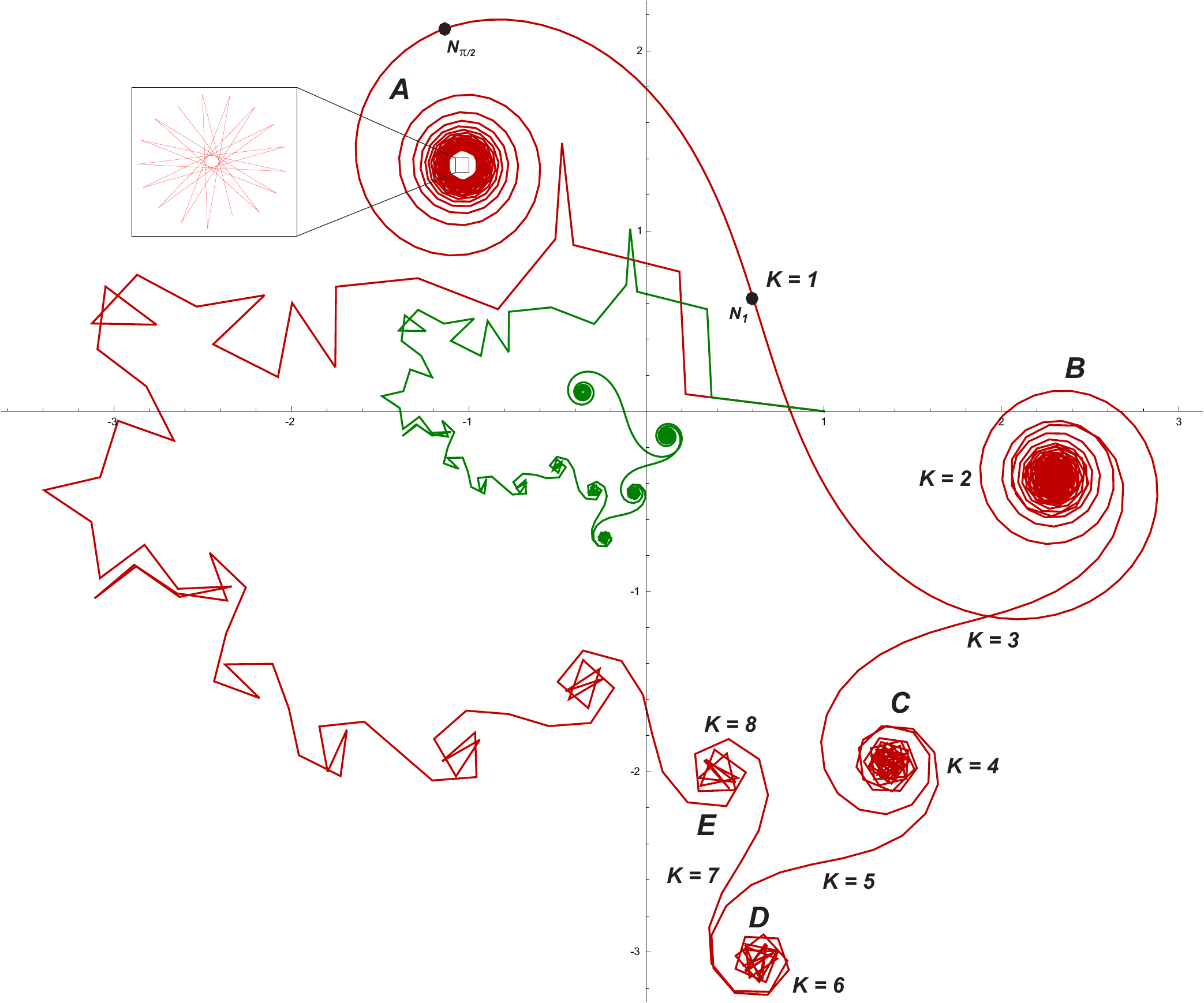}
	\caption{Paths described by the partial sums of the series for $\eta(s)$\, and $\eta(1-\overline{s})$\,, $\; s=\frac{1}{2}-\alpha+it$\,.}
	\label{fig:1}       
\end{figure}

Imagine now to walk along one such path in the direction of increasing $N$. Upon arriving at the end of segment $u_N(s)$,  in order to step on segment $u_{N+1}(s)$ we will then need to turn left or right, according to the value of the difference between their respective  arguments, noting that $\theta_n(t) = -t \log n $  when $n$ is odd, whereas  $\theta_n(t) =\pi -t \log n $ when $n$ is even. Thus, the included angle between $u_N(s)$ and $u_{N+1}(s)$ is $\delta_{N+1}(t) = t \log \frac{N+1}{N}$. Suppose we have already walked our way through that kind of "whirl" labeled as D, and we are thus turning to our right when stepping on the next segment, as soon as we arrive at the end of a segment for which $\delta_{N+1}(t)$ switches from the range of values $K \pi>\delta_{N+1}(t)> (K-1) \pi$, $K=6$, to the range of values $5\pi>\delta_{N+1}(t)> 4\pi$ (i.e. $K=5$), we will need to turn to the left. We would then return to take turns to the right somewhere in the middle of whirl C, when $\delta_{N+1}(t)$ crosses from range $K = 5$ into range $4\pi>\delta_{N+1}(t)> 3\pi$ (i.e. $K=4$), and so on and so forth, until finally reaching the $K=1$ range $\pi>\delta_{N+1}(t)> 0$, then forever taking turns to the left, while getting closer and closer to $\eta(s)$ (located somewhere in the middle of whirl A, which we can hence name \textit{End Whirl}). The square inset enlarges a detail (for just 16 segments) as it would appear by zooming deep inside the \textit{End Whirl}, revealing how it eventually settles into a very simply structured star-shaped pattern. K-range crossing from $K$ into $K-1$ occurs at the end of the $u_{N_K}$ segment, where
\begin{eqnarray}
	N_K = \lceil \frac{1}{e^{\frac{K \pi}{t}}-1} \rceil
	\label{eq_7}
\end{eqnarray}
Incidentally, the study of the size of said K-ranges may suggest a sufficient condition for the Lindelöf Hypothesis. The number of K-ranges is $K_{MAX} = t \log 2 / \pi$\, (inverting \ref{eq_7} when $N_K = 1$). Were it possible to show that along the critical line the size of each K-range is bounded by $\frac{C}{k}$\, , $C>0$ no matter how large, then:   
\begin{equation}
	\left| \eta(s) \right| < C \sum_{k=1}^{K_{MAX}}\frac{1}{k} = C \left( \log K_{MAX} + \gamma + ..... \right) = C \log t + C \gamma + C \log \frac{\log 2}{\pi} + ..... \; \rightarrow \; O\left( t^\epsilon \right) 
	\label{eq_8}
\end{equation}
Note that such $O\left( t^\epsilon \right)$\, behavior (Lindelöf Hypothesis for the Riemann Zeta function, but which directly transfers to the $\eta$ function thanks to the fact that $ \sqrt{3-2\sqrt{2}}  \leq |1-2^{0.5-it}| \leq \sqrt{3+2\sqrt{2}}$ ) would still be preserved even if the constant $C$\, were instead a function of $t$ growing no faster than $C(t) = c \left( \log t \right) ^{b} $, $b>0$ and $c>0$ no matter how large. 

\section{The Geometric Approach}
\label{sec:3}
By studying the geometrical pattern traced by the $\eta_n(s)$, it is possible to identify stricter bounds for the asymptotic size of the remainder $R_n(s)$. Said asymptotic size is already known to be of the order of its first term (for example through the application of Euler MacLaurin summation formula). Here we will instead follow a geometric approach resulting in a strict asymptotic relationship, while also producing some other interesting insights. \medskip \\
\textbf{Theorem 1.\;}\itshape For $s=\sigma + it$\, in the interior of the critical strip, denote by  $ \eta_{n-1}(s)$\, the $(n-1)^{th}$\, partial sum of the series for the $\eta$\, function, by  $R_n(s)$\, the corresponding remainder, and by $u_n(s)$\, the line segment starting at $\eta_{n-1}(s)$\, to end at $\eta_n(s)$\,. For sufficiently large $n$\,, it follows that 
\begin{equation}
	the\; circle\; of\; diameter\; u_{n+2}(s)\; lies\; strictly\; inside\; that\; of\; diameter\; u_{n}(s)\;\; , \nonumber 
\end{equation}
hence 
\begin{equation}
	|R_n(s)|\;\; < \;\; \frac{1}{n^{\sigma}} \; \label{eq_9}
\end{equation}
moreover, a strict asymptotic relationship holds 
\begin{equation}
	as\;\;\;\; n \rightarrow \infty \;\;\;\;\;\;\;\;  R_n(s)\;\; = \;\; \frac{(-1)^{n-1}}{2\,n^{s}}\; + \; \epsilon_n(s) \;\;,\;\;\;\; \frac{\epsilon_n(s)}{n^{-\sigma}} \;\; \rightarrow \;\;\ 0    \label{eq_10}
\end{equation}
actually, one even has $\epsilon_n(s) \; \in \; o\left( \frac{1}{n^2}\right)$.
\medskip \\
\textbf{PROOF. \;}\normalfont  

The procedure followed for this proof is based on a geometric analysis of the pattern of convergence displayed by the path traced by the partial sums $\eta_n(s)$\,. Said path is composed of line segments, and for simplicity we will use the term $u_n$\, to refer to the line segment starting at $\eta_{n-1}(s)$\, and ending at $\eta_n(s)$\,, while implicitly assuming its dependence on $s$. Furthermore, because pair of paths corresponding to  $s=\sigma+it$\, and $\overline{s}=\sigma-it$\, are mirror symmetrical with respect to the real axis, it will be sufficient to prove the following results for the case $t > 0$\,.\\ 
\indent
While approaching their respective point of convergence, $\eta(s)$\,, paths corresponding to different values of $\sigma > 0$\, appear all to end up describing simply structured star-shaped patterns (recall the inset in Fig. 1), characterized by angles between consecutive segments being $<\frac{\pi}{2}$\,. This is actually the result of having to add $\pi $\,, to the argument $-t\ln(n) $\,, every other segment (because of the alternating sign). In fact, when $ n $\, becomes sufficiently large, $t\ln(n+1) $\, will be just a bit larger than $t\ln(n) $\,, and because one of the two segments will need to be turned around by $\pi $\, (the segment corresponding to even $ n $\,), the angle between said consecutive segments will eventually become an acute angle \, $\delta_{n+1} = t\ln(n+1) - t\ln(n) < \frac{\pi}{2}$\, (Fig. 3a helps in visualizing this), shrinking down more and more as $ n $\, grows larger and larger. Being interested solely in the absolute value of said acute angle, it can be easily verified that its value is $|-t\ln\frac{n+1}{n}|=|t\ln\frac{n}{n+1}|$\,, which is the same as $t\ln\frac{n+1}{n}$\, when $t \geq 0$\,, and hence $\rightarrow 0$\, as $n\rightarrow\infty$\,. Ultimately, the observed star-shaped pattern is hence a direct result of the alternating signs in (\ref{eq_2}). \\
\begin{figure}
	\centering
	\subfigure[The circles of radius $r_n$\,, $r_{n+2}$\,, $\frac{r_n}{2}$\,, $\frac{r_{n+2}}{2}$\,.]
	{\includegraphics[width=0.4\textwidth]{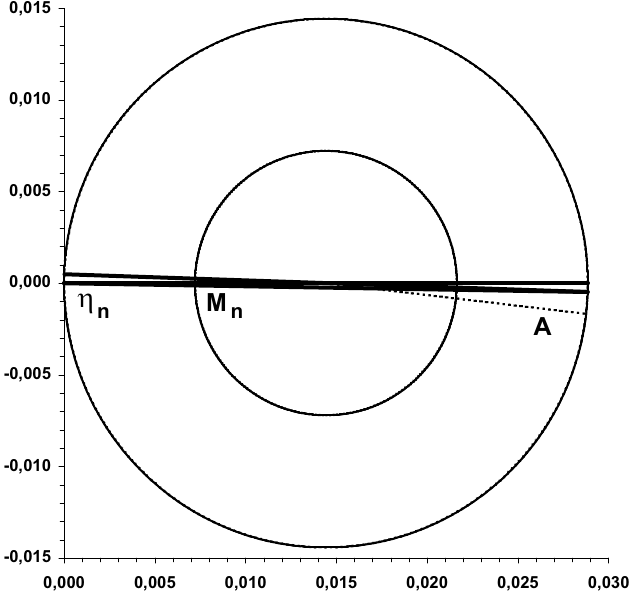}}
	\hspace{15mm}
	\subfigure[Expanding tenfold the vertical scale.]
	{\includegraphics[width=0.4\textwidth]{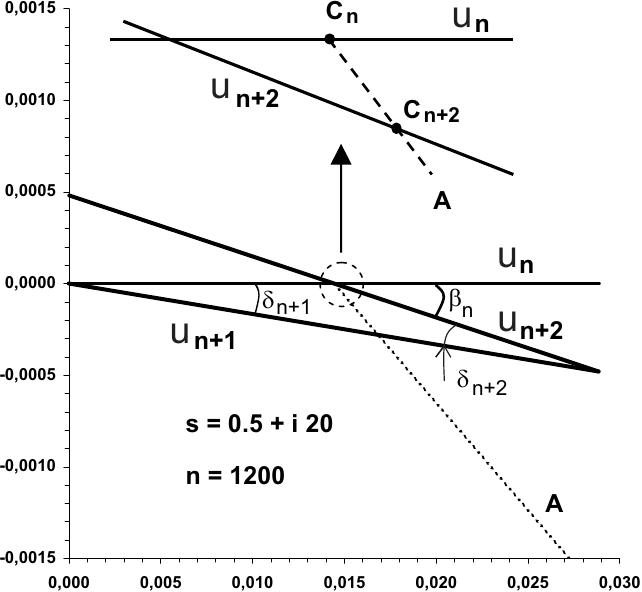}}
	\subfigure[Disk of radius $r_{n+2}$\, $\subset$\, disk of radius $r_n$\,.]
	{\includegraphics[width=0.4\textwidth]{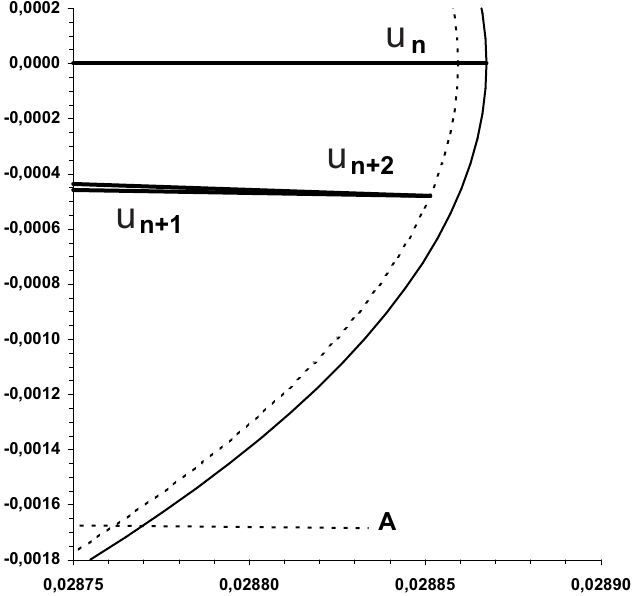}}
	\hspace{15mm}
	\subfigure[Disk of radius $\frac{r_{n+2}}{2}$\, $\subset$\, disk of radius $\frac{r_n}{2}$\,.]
	{\includegraphics[width=0.4\textwidth]{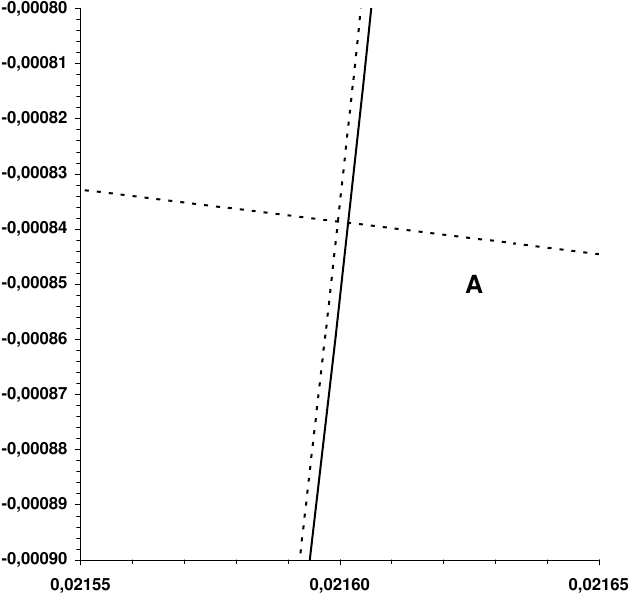}}
	\caption{proving that $|R_n(\sigma + i t)|\, \sim \, \frac{1}{2n^{\sigma}}$\, as $n\rightarrow\infty$\,.}
	\label{fig:2}
\end{figure} 
\indent
Fig. 2 illustrates an example, limited to just 3 consecutive segments, of how said pattern appears once the path described by the partial sums has already settled into the star-shaped pattern at the center of the \textit{End Whirl} (Fig. 1). We will at first demonstrate that said star-shaped pattern ought to be bounded in size. 
To this end we construct a circle whose diameter is $u_n$\,, and whose center, $C_n$\,, is its the midpoint. Then, if the disk defined by the circle constructed in a similar way on $u_{n+2}$\, were eventually (i.e.: as $n\rightarrow\infty$\,) contained within the disk defined by the circle constructed on $u_n$\,, we would have proved that the distance between point $\eta_n(s)$\, and the point of convergence, $\eta(s)$\,, cannot be greater than the length of $u_{n}$\,. Fig. 2a illustrates the geometric construction on which this proof will be based. As this proof is concerned solely with the geometric relationships  of $u_n$\,, $u_{n+1}$\,, and $u_{n+2}$\,, relatively to each other, our frame of reference can be chosen in any convenient way which preserves said relative properties. The origin of the frame of reference chosen for Fig.2 coincides with the point $\eta_n(s)$\,, while the real axis is chosen as to overlap $u_{n}$\,. In Fig. 2b the scale of the imaginary axis has been purposely stretched by a factor of 10, only to make it easier to visualize the pattern described by the partial sums. Referring to the triangle defined by the intersections of $u_{n}$\,, $u_{n+1}$\,, and $u_{n+2}$\,, let us denote by $\delta_{n+1}=t\ln\frac{n+1}{n}$\, the angle between $u_{n+1}$\, and $u_{n}$\,, by $\delta_{n+2}=t\ln\frac{n+2}{n+1}$\, the angle between $u_{n+2}$\, and $u_{n+1}$\,,, and by $\beta_n$\, the supplement of the third angle of said triangle (i.e. the angle between $u_{n}$\, and $u_{n+2}$\, ). By construction, angle $\beta_n$\, corresponds to $\beta_n=\delta_{n+1}+\delta_{n+2}= t \ln\frac{n+2}{n}$\,.\\
\indent

Fig. 2b shall then be read as follows: \\[0.1in]
$u_{n}$\, starts off at $\frac{1}{n^{\sigma}}+i0$\,, ending at $0+i0$\, ;  \\[0.1in]
$u_{n+1}$\, starts off at $0+i0$\,, ending at $\frac{1}{(n+1)^{\sigma}}\left(\cos \delta_{n+1}-i\sin\delta_{n+1}\right)$\,;  \\[0.1in]
$u_{n+2}$\, starts from the end of \textit{segment} $n+1$\,,  ending at $\frac{\cos\delta_{n+1}}{(n+1)^{\sigma}}-\frac{\cos\beta_n}{(n+2)^{\sigma}}+i\left(-\frac{\sin\delta_{n+1}}{(n+1)^{\sigma}}+\frac{\sin\beta_n}{(n+2)^{\sigma}}\right)$\,. \\[0.1in]   

Line \textbf{A} is drawn through $C_n$\, and $C_{n+2}$\,, the midpoints of $u_{n}$\, and $u_{n+2}$\, (see the enlarged detail in Fig. 2b). $C_n$\, is also the center of the circle (solid line) having $u_{n+2}$\, as its diameter, and hence the radius of said circle is $r_n=\frac{1}{2n^{\sigma}}$\,. $C_{n+2}$\, is also the center of the circle (dotted line) having $u_{n+2}$\, as its diameter, and hence the radius of said circle is $r_{n+2}=\frac{1}{2(n+2)^{\sigma}}$\,. The purpose of line \textbf{A} is to make it easier to identify the points of closest approach of the two circles. In Fig. 2c said circles appear stretched. This is the result of having chosen different scale factors for the two axes, so that line \textbf{A} may also appear in the graph together with all other segments. It is clear that for the example of Fig.2 ($s=\frac{1}{2} + i 20$\,, $n=1200$\,) the disk of radius $r_{1202}$\, is contained within the disk of radius $r_{1200}$\, (closest approach is their distance along line \textbf{A}). Several other examples were verified by numerical simulations, always resulting in \textit{disk} $n+2$\, being eventually contained strictly inside \textit{disk} $n$\,, and therefore suggesting that this kind of asymptotic behavior might indeed be representative of the most general case. To demonstrate that this is actually the case, any effective proof will need to rely solely on assumptions eventually satisfied, as $ n\rightarrow\infty$\,, by the partial sums $\eta_n(\sigma+it)$\,. As a matter of fact, for values $n$\, such that $t\ln\frac{n+1}{n}<\frac{\pi}{2}$\,, any three consecutive segments will define a triangle of the same kind as the one depicted in Fig. 2b, and this irrespective of the particular value $\sigma>0$\, (more precisely, the triangle in Fig. 2b actually refers to  $\delta_{n+1}$\, and $\delta_{n+2}$\, both $<\frac{\pi}{4}$\,, but it is easy to verify that (\ref{eq_11}), identifying the locations of $C_n$\, and $C_{n+2}$\,, are valid even when said angles are both just a tiny little bit less than $\frac{\pi}{2}$\,). We shall also remark that the validity of the following proof cannot possibly be affected by our choice to start from a line segment of even index, rather than from one of odd index. \\

The following proof will refer solely to values $n\geq N$\,, where $N=N(t)$\, is the smallest integer that is not less than $1 / \left(e^{\frac{\pi}{2t}}-1\right)$\,. Let us first evaluate the distance $\Delta_c = \Delta_c(n) = \left|C_{n+2}-C_{n}\right|$  
\begin{equation}
	C_n = \frac{1}{2n^{\sigma}}+i0  \; \; \; \; \; \; \; \; \; \; \; \;  C_{n+2} = \frac{\cos\delta_{n+1}}{(n+1)^{\sigma}}-\frac{\cos\beta_n}{2(n+2)^{\sigma}}+i\left(-\frac{\sin\delta_{n+1}}{(n+1)^{\sigma}}+\frac{\sin\beta_n}{2(n+2)^{\sigma}}\right) \nonumber
\end{equation}
\begin{eqnarray}
	&   &\Delta_c^{2} =  \left|C_{n+2}-C_{n}\right|^{2} =   \nonumber \\ 
	& = & \frac{\cos^{2}\delta_{n+1}}{(n+1)^{2\sigma}} + \frac{\cos^{2}\beta_n}{4(n+2)^{2\sigma}} - \frac{\cos\delta_{n+1} \; cos\beta_n}{(n+1)^{\sigma}(n+2)^{\sigma}} + \frac{1}{4n^{2\sigma}} + \frac{\cos\beta_n}{2n^{\sigma}(n+2)^{\sigma}} - \frac{\cos\delta_{n+1}}{n^{\sigma}(n+1)^{\sigma}} \nonumber \\
	& + & \frac{\sin^{2}\delta_{n+1}}{(n+1)^{2\sigma}} + \frac{\sin^{2}\beta_n}{4(n+2)^{2\sigma}} - \frac{\sin\delta_{n+1} \; \sin\beta_n}{(n+1)^{\sigma}(n+2)^{\sigma}} \nonumber \\
	& = & \frac{1}{(n+1)^{2\sigma}} + \frac{1}{4(n+2)^{2\sigma}} - \frac{\cos(\beta_n-\delta_{n+1})} {(n+1)^{\sigma}(n+2)^{\sigma}} + \frac{1}{4n^{2\sigma}} + \frac{\cos\beta_n}{2n^{\sigma}(n+2)^{\sigma}} - \frac{\cos\delta_{n+1}}{n^{\sigma}(n+1)^{\sigma}} \;\; . 
	\label{eq_11}
\end{eqnarray}
For the difference between the two radiuses, $\Delta_r = \Delta_r(n) = \left|r_{n}-r_{n+2}\right|$\,, we have  
\begin{eqnarray}
	\Delta_r^{2} & = & \frac{1}{4n^{2\sigma}} + \frac{1}{4(n+2)^{2\sigma}} - \frac{1}{2n^{\sigma}(n+2)^{\sigma}} \; . \label{eq_12}
\end{eqnarray}
Subtracting (\ref{eq_11}) from (\ref{eq_12}) and simplifying, $\Delta_r^{2} - \Delta_c^{2}$\, becomes
\\
\begin{equation}
	\frac{\left(1+\frac{1}{n}\right)^{\sigma}\left[\left(n+2\right)^{\sigma}cos\left(\delta_{n+1}\right) - \left(n+1\right)^{\sigma}\frac{1+cos\left(\beta_n\right)}{2}\right] - \left[\left(n+2\right)^{\sigma} - \left(n+1\right)^{\sigma} cos\left(\beta_n - \delta_{n+1}\right)\right]}
	{\left(n+1\right)^{2\sigma}\left(n+2\right)^{\sigma}} .
	\label{eq_13}
\end{equation}

\indent
\\
Considering that for any given $t$\, the cosine factors  $\rightarrow 1 \, + \, O(1/n^2)$\, as $n\rightarrow\infty$\,, both the square brackets at the numerator expand to $\left[ (n+2)^\sigma - (n+1)^\sigma +\, O(n^{1-2/\sigma}) \right] $\,, valid for any $0<\sigma<1$. Thus, the numerator is eventually a strictly positive quantity. However, at smaller $n$\, values the oscillatory nature of those cosine factors has a much greater impact, and it is not difficult to verify the existence of values of $n$\, resulting in (\ref{eq_13}) $<\,0$\,. Thus, because (\ref{eq_13}) $\rightarrow\ 0$\, from eventually positive values, it means that there must exist a value of $h = h(\sigma, t)$\, at which (\ref{eq_13}) definitely crosses over from negative to positive values.
\begin{figure}
	\centering
	\subfigure[For $n$\, sufficiently large,  $\delta_n<\frac{\pi}{2}$\,.]
	{\includegraphics[width=0.4\textwidth]{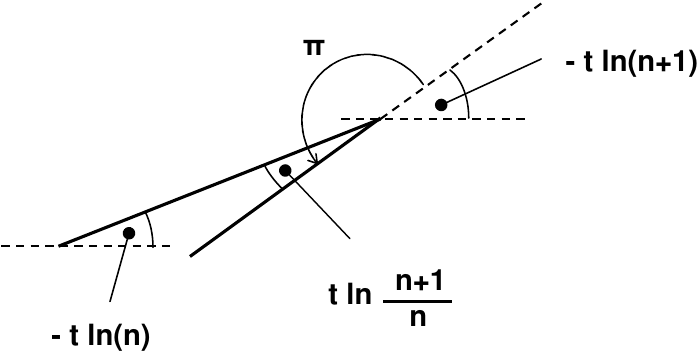}}
	\hspace{15mm}
	\subfigure[plotting (\ref{eq_13}) for $t=300 \, , \, \sigma=0.65$\,.]
	{\includegraphics[width=0.425\textwidth]{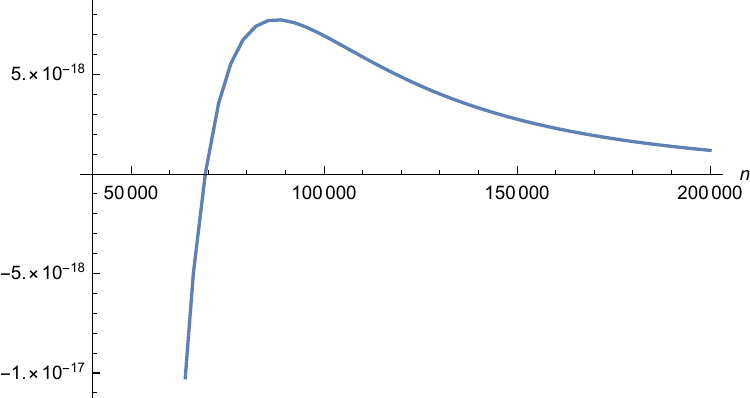}}
	\caption{}
	\label{fig:3}
\end{figure}
Fig. 3b plots (\ref{eq_13}) for the example $t=300 \, , \, \sigma=0.65$\,, while for the general case it can be verified that said last crossing from negative to positive values occurs around $\frac{t^2}{2\sigma} $\, (though tedious, it can be verified by expanding at $x=\frac{t^2}{2\sigma} $\, to then observe the change of sign when increasing from $n=\lfloor x \rfloor - 1$\, to $\lceil x \rceil + 1$ ), although such result will not be used here. For all values $n>h$\, it is therefore $\Delta_r^{2} - \Delta_c^{2}> 0$\,, implying that $\Delta_r$\, is eventually $\Delta_r > \Delta_c$\,, in turn implying that the circle built on $segment \; n+2$\, is strictly contained within the circle built on $segment \; n+2$\,, the one built on $segment \; n+4$\, is strictly contained within that built on $segment \; n+2$\,, and so on, as $n\rightarrow\infty$\,, while shrinking down to the point of convergence, $\eta(s)$\,, and thus proving (\ref{eq_9}). \\
By tediously expanding (\ref{eq_11}) for $n\rightarrow\infty$\, it turns out that terms of the order $n^{-2\sigma}$ all cancel each other out, so that the next remaining terms are those of the order $n^{-4-2\sigma}$. Hence, 
\begin{equation}
	\left|C_{n+2}-C_{n}\right|\,=\, O \left( \frac{1}{n^{2+\sigma}}\,\right) \label{eq_14}
\end{equation}
whereas, by definition, the diameters of the nested circles vanish at the $n^{-\sigma}$ rate. Thus, the distance (\ref{eq_14})  between the centers of two consecutive nested circles vanishes at a much faster rate than their diameters, in turn implying that the nested circles ought to be asymptotically concentric, finally proving (\ref{eq_10}). \medskip \qed \\
In fact, by observing that $R_n(s)$ vanishes closer and closer to both $C_{n}$ and $C_{n+2}$, implying that the error term $R_n(s) - n^{-s}/2$ cannot shrink slower than (\ref{eq_14}), by further remarking that $\sigma > 0$ we could as well conclude that $|R_n(s) - n^{-s}/2| \in o(n^{-2})$.  \medskip \\ 
To aid in visualizing the result of \textbf{Theorem 1}, Fig. 4a illustrates the star shaped pattern drawn by connecting up consecutive values of $\eta_k(s_0)$\ (i.e. the vertexes), for $s_0 = 1/2\, +\, i \,\, 37.5861781588256$\,, the sixth non-trivial zero (at zeros, it is of course $\eta_{k-1}(s_0)=-R_{k}(s_0)$). It depicts solely partial sums featuring $k\,>\,\frac{t_0^2}{2\sigma_0} $\,, so that the circle of diameter $u_{k+2}(s_0)$\, is strictly contained within the circle of diameters $u_{k}(s_0)$\,. It is apparent how $\eta(s_0)=0$\, is indeed located somewhere very near the centers of said nested circles. The pattern of convergence towards the "center" (using "center" to roughly indicate the point of convergence) of the star-shaped pattern can then be best appreciated when plotting only even-indexed partial sums (the grey dots in Fig. 4b), rather than each one of them (which will jump back and forth between nearly diametrically opposed positions, because of the alternating sign, as clearly illustrated by Fig. 4a). Similarly, this also holds when plotting only the odd-indexed ones.
\medskip \\
Thus, without loss of generality, we will henceforth restrict our discussion to even values of $N$.
\medskip \\
The smaller black dots in Fig. 4b represent instead the values computed using (\ref{eq_10}). By zooming in on Fig. 4b (the figures in this PDF allow zooming to a very fine level of detail), it becomes visually apparent that, as one moves toward the inner turns of the spiraling pattern (i.e., as $k$ increases), the black dots overlap the gray ones with increasing accuracy. This provides numerical evidence supporting the result of \textbf{Theorem 1}.
For convenience, in the following we will use the term \textit{Deep Spiral} $N$ to denote the pattern of convergence of the $R_{2k}(s)$, starting from index $2k = N$, where $N$ is any sufficiently large value located deep inside the \textit{End Whirl} and satisfying the condition "the circle of diameter $u_{N+2}(s)$\; lies strictly inside that of diameter $u_{N}(s)\;$". Upon closer inspection, Fig. 4b reveals that the spacing between dots decreases as $k$ increases, both in the radial direction and along the spiral itself. In polar coordinates, this corresponds to an increase in the density of dots along both the radial and angular coordinates. It is not difficult to show that this is a general result, independent of the particular choice of $s_0$ used in Fig. 4b. Let us consider the first full turn of \textit{Deep Spiral} $N$, for sufficiently large $N$, the asymptotic expression (\ref{eq_10}) suggests that one $2\pi$ revolution along the spiral would very nearly correspond to a $2\pi$ rotation of $-1/(2k)^{s}$. Thus, starting at $N$, said $2\pi$ angular span is covered by index:
\begin{equation}
	k = N+2j \;, \;\;\;\; 0\leq j \leq m \;, \;\;\;\; m = \lfloor N \left( e^{\frac{2 \pi}{t}}-1 \right)/2 \rfloor \; \xrightarrow[t \rightarrow \infty]{} \; N\, \frac{\pi}{t} \;, \;\;\;\; N > \frac{t^2}{2 \sigma}  \label{eq_15}
\end{equation}
thus, $m$ corresponds to the number of dots present in the first full $2\pi$ turn of \textit{Deep Spiral} $N$ (more exactly, the angular span nearest to $2\pi$ and $\leq 2\pi$). Notice how, at large $t$ values, the length of the first turn of \textit{Deep Spiral} $N$ contracts at nearly the same rate as the circumference of the circle of diameter $1/N^{\sigma}$. We can hence estimate that the density of dots along the spiraling pattern grows at nearly the same rate as  $N^{1+\sigma} / t$ (while it is $N > \frac{t^2}{2 \sigma}$). Albeit at a slower pace, also the radial density increases, given that the radial distance between two consecutive turns would contract as $2 \, \Delta_r \simeq 1/N^{\sigma}-1/(N+2m)^{\sigma} \; = \; \frac{1}{N^{\sigma}} \left( 1 - \left( 1 + \frac{2m}{N}\right)^{-\sigma} \right)$, which, at large $t$ values, is well approximated by $\Delta_r \simeq \pi \, \sigma /(t\,N^{\sigma})$. As an example, for $t>1000$ and along $\sigma = 1/2$, that would correspond to $\Delta_r < $ about 0.3\% of $r$. Hence, at sufficiently large $t$, the radial density is well approximated by $1/\Delta_r \simeq N^{\sigma} \frac{t}{\pi\sigma}$ (while it is $N > \frac{t^2}{2 \sigma}$).  
\medskip \\ 
We shall now remark how the remainder $R'_k(s)$ of the series for the first derivative of the Dirichlet $\eta$ function,  $\eta'(s)=\eta'_{k-1}(s) + R'_k(s)$, also follows a spiraling pattern very similar to that in Fig 4b. Indeed: \medskip          
\begin{figure}
	\centering
	\subfigure[the $\eta_{k}(s_0) \longrightarrow \eta(s_0) = 0$]
	{\includegraphics[width=0.48\textwidth]{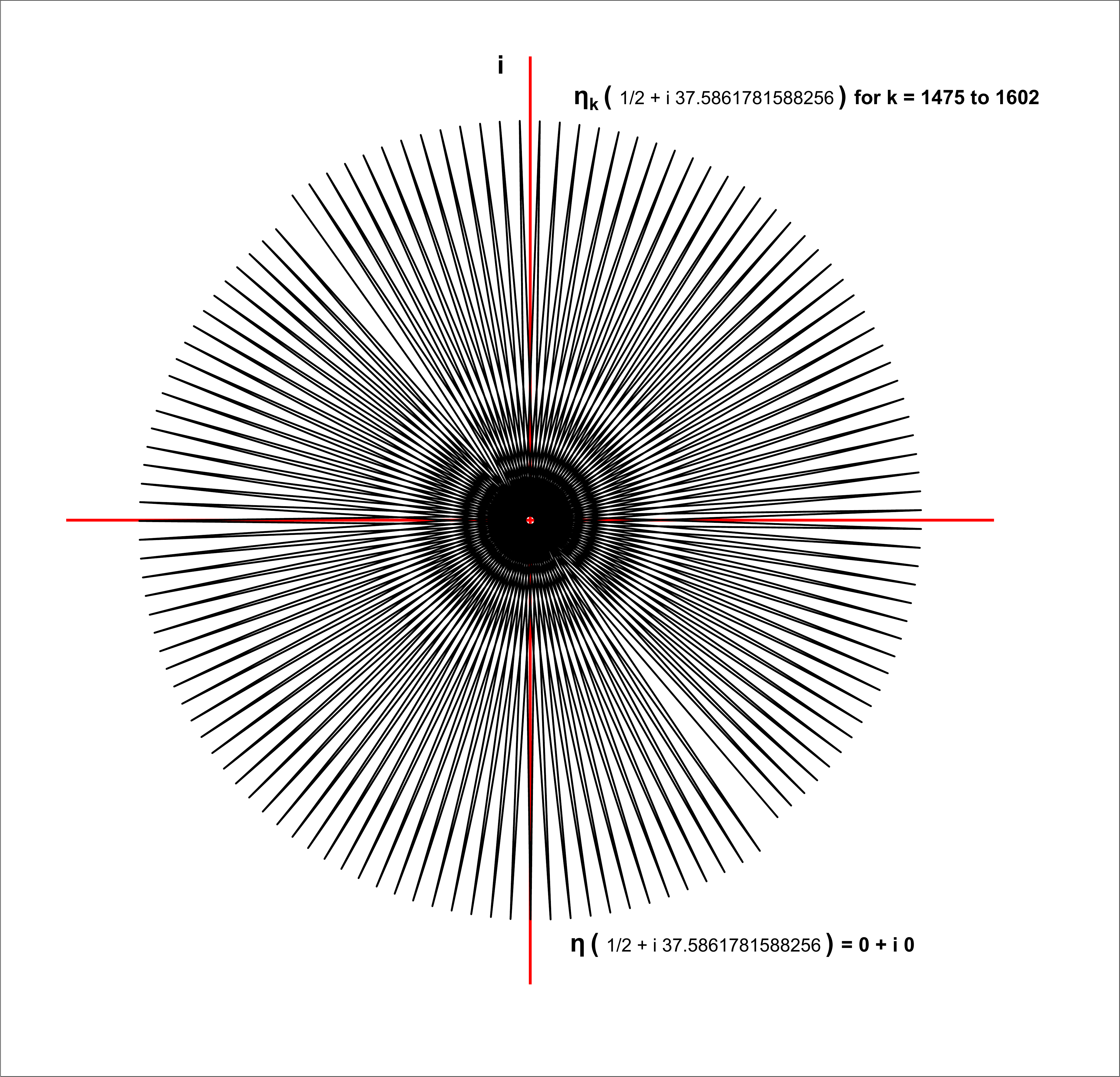}}
	\subfigure[even-indexed $\eta_{2k-1}(s_0)=-R_{2k}(s_0)$ for \textit{Deep Spiral 1476}]
	{\includegraphics[width=0.48\textwidth]{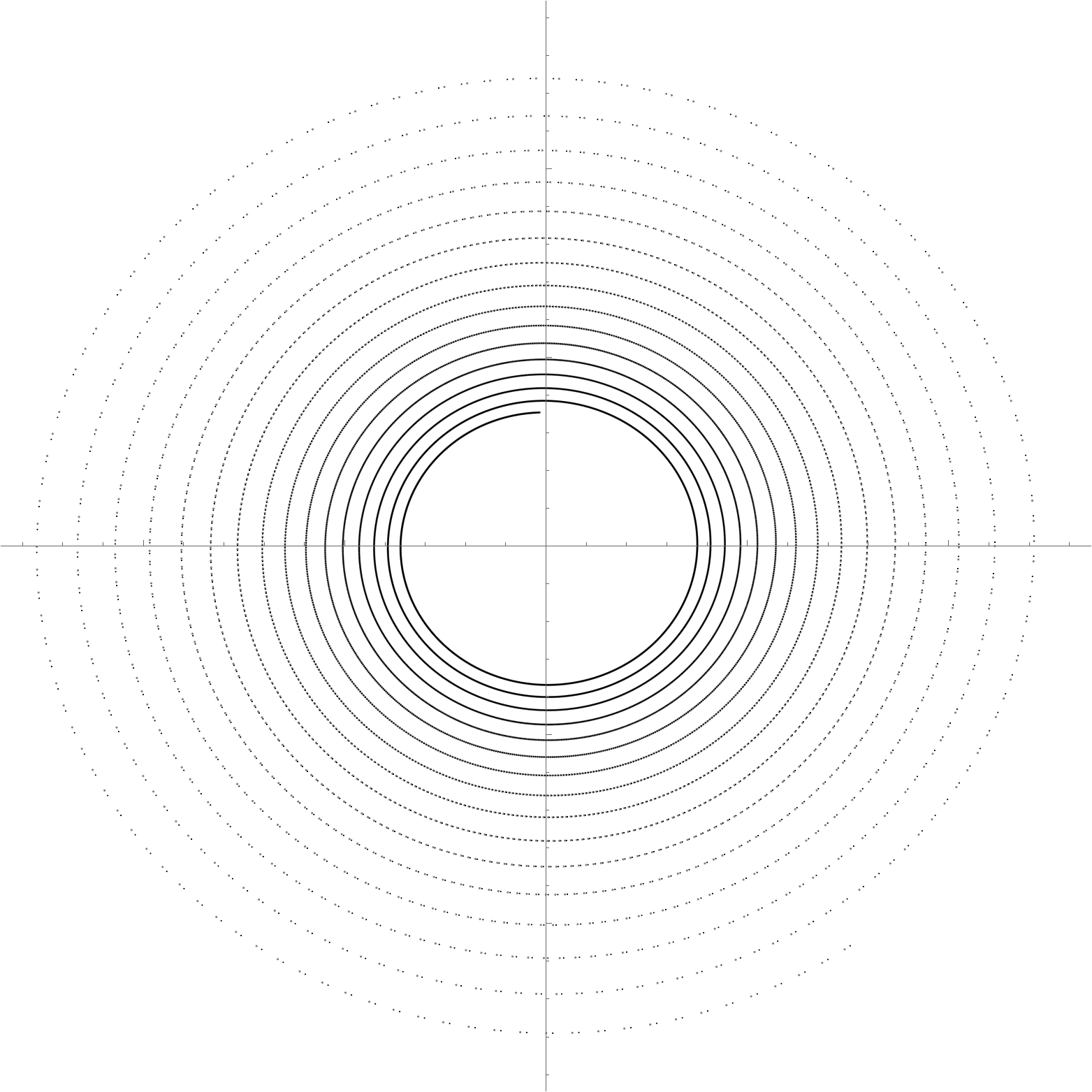}}
	\caption{$\eta_k(s)$ examples at a known zero. At the right, when joining only even indexed $\eta_k(s)$.}
	\label{fig:4}
\end{figure}

\textbf{Corollary 1.\;}\itshape For $s=\sigma + it$\, in the interior of the critical strip, denote by  $ \eta'_{n-1}(s)$\, the $(n-1)^{th}$\, partial sum of the series for the first derivative of the $\eta$\, function, and by  $R'_n(s)$\, the corresponding remainder. The following strict asymptotic relationship holds 
\begin{equation}
	as\;\;\;\; n \rightarrow \infty \;\;\;\;\;\;\;\;  R'_n(s)\;\; = \;\; \frac{(-1)^{n} \log n}{2\,n^{s}}\; + \; \delta_n(s) \;,\;\;\; \frac{\delta_n(s)}{n^{-\sigma} \log n} \; \rightarrow \; 0   \label{eq_16}
\end{equation}
actually, one even has $\delta_n(s) \; \in \; o\left( \frac{\log n}{n^{2}}\right)$.
\medskip \\
\medskip \\ 
\textbf{PROOF. \;}\normalfont As $n \rightarrow \infty$, it is $\log (n+1) = \log n + 1/n + \cdots$ and $\log (n+2) = \log n + 2/n + \cdots$. Retracing then the same geometric proof as for \textbf{Theorem 1}, we will obtain, for both the $\Delta_C$ and the $\Delta_r$ expressions, a first part composed of all terms multiplied by the common factor $\log n$, added to a second part where the same terms are instead multiplied by either $1/n$ or $2/n$, and which will hence vanish at a much faster rate than said first part.  \qed 
\medskip \\
Because at sufficiently large values of $t$ and $N$ we have $\log(N+2m)/\log N \, \sim \, 1 + \log(1+\frac{2 \pi}{t})/\log N$ (see (\ref{eq_15}) above), for the first turn of \textit{Deep Spiral} $N$ it will be $\log(N+2m) \simeq \log N$, so we may just as well replace $\log n$ in (\ref{eq_16}) by $\log N$ (or, for that matter, by $\log(N+2m)$) and still maintain good accuracy. Indeed, if we were to plot the \textit{Deep Spiral} $N$ for the $ R'_n(s)/\log N$, at a sufficiently large $N$ we will observe an almost perfect overlap with the \textit{Deep Spiral} $N$ for the $ R_n(s)$. 
\medskip \\         
Moreover, if $\eta(s_0)=0$, and $\Re(s_0)>0 $, observing that the sequence of partial sums $\lbrace \eta_{k}(z)\rbrace$ converges uniformly on compact subsets of the half-plane $\Re(z)  > 0$ to $\eta(z)$, and that for $\eta(z)$ it is possible to find a sufficiently small circular boundary $|z-s_0|= \epsilon$ surrounding only the $s_0$ zero (zeros are isolated), while along said boundary it is also $\eta(z) \neq 0$, a known result due to Hurwitz \cite{06} implies that there is an integer $N$ such that for $n\geq N$, $\eta(z)$ and $\eta_n(z)$ have the same number of zeros inside the domain encircled by $|z-s_0|= \epsilon$ (counting multiplicities is implied). There must hence exist a sequence $\lbrace z_k(s_0)\rbrace$ such that $\eta_{k-1}(z_k(s_0)) = 0$, whereby $z_k(s_0) \rightarrow s_0$\, as $k\rightarrow\infty$, $s_0$ representing an accumulation point for the sequence $\lbrace z_k(s_0)\rbrace$. Denoting by $s_0 = \sigma_0 + i t_0$ a zero of the $\eta$ function, and referring to the even-indexed Hurwitz zeros sub-sequence $\lbrace z_{2k}(s_0)\rbrace$, $z_{2k}(s_0) \rightarrow s_0$ as $k \rightarrow \infty$, defined so that $\eta_{2k-1}(z_{2k}) = 0$, and which in general does not correspond to $\eta_{2k-1}(1 - z_{2k}) = 0$, by definition implies that $\eta(z_{2k}) = R_{2k}(z_{2k})$ (see \ref{eq_4}). Now, trying to directly prove that sequences $\lbrace z_{2k}(s_0)\rbrace$, featuring the properties just described, can only exist when $\sigma_0 = 0.5$ appears challenging. However, by combining the results of \textbf{Theorem 1} with the constraints impose by the Functional Equation yields a perhaps easier to approach strict equality.    \medskip \\
\textbf{Corollary 2.\;}\itshape In the interior of the critical strip, the falsehood of the Riemann Hypothesis requires the existence of at least one sequence $\lbrace z_{2k}(s_0) \rbrace$, $\Re(s_0) \neq 1/2$, such that the following equality is eventually satisfied
\begin{equation}
	\eta_{2k-1}(1-z_{2k}(s_0)) = \chi^{\pm}(z_{2k}(s_0)) R_{2k}(z_{2k}(s_0)) - R_{2k}(1-z_{2k}(s_0)) \label{eq_17}
\end{equation}  
\medskip \\
\textbf{PROOF. \;}\normalfont
\\
By the definition of $\lbrace z_{2k}(s_0)\rbrace$, it is $R_{2k}(z_{2k})= \eta(z_{2k})$. The functional equation $	\eta(1-z_{2k}) = \chi^{\pm}(z_{2k}) \eta(z_{2k})$, considering that $\eta(1-z_{2k}) =	\eta_{2k-1}(1-z_{2k}) + R_{2k}(1-z_{2k})$, leads then to (\ref{eq_17}) in a straightforward way.\qed
\medskip \\   
Conversely, if one can prove that no sequence of Hurwitz zeros ${z_{2k}(s_0)}$ satisfies the strict equality above, then the Riemann Hypothesis is confirmed. To that aim, a sufficient condition would be to show that the asymptotic behaviours of the two sides of (\ref{eq_17}) cannot possibly match as $k \to \infty$.\\
Denote $z_{2k}(s_0) = s_0 + \delta_{2k}(s_0)$. To rewrite the strict equality above in terms of the asymptotic behaviour of the vanishing left- and right-hand terms, we use the power series expansion of the analytic function on the left-hand side together with the result of \textbf{Theorem 1}. This yields the following asymptotic equivalence (the dots denote higher-order terms vanishing at faster rates):
\begin{equation}
	\eta_{2k-1}(1-s_0)-\eta'_{2k-1}(1-s_0) \delta_{2k} + \cdots = - \frac{\chi^{\pm}(s_0+\delta_{2k})}{2\,(2k)^{s_0+\delta_{2k}}} + \frac{1}{2\,(2k)^{1-s_0-\delta_{2k}}} + \cdots \label{eq_18}
\end{equation}  
At present, however, this seems to be little more than an interesting curiosity, as it remains unclear whether it can meaningfully aid in studying the Riemann Hypothesis. Any progress in this direction would first require a detailed analysis of the asymptotic behaviour of the two partial sums on the left, a daunting task indeed.\\     

\section{A Reformulation of The Riemann Hypothesis}
\label{sec:3}

\textbf{Theorem 2. } \itshape The statement that $\lim_{n\to\infty} X^{\pm}_n(s)=L(s)$\, is continuous on $D=\left\lbrace s \in \mathbb{C}: \; 0< \Re(s) < \frac{1}{2}\right\rbrace $\, is a necessary and sufficient condition for the following inequality to be satisfied in $D$\,
\begin{equation}
	\eta(s) \neq 0. \nonumber
\end{equation} 
\\ 
\textbf{PROOF.} \normalfont
In order to prove necessity and sufficiency, we are going to follow these two logical steps:
\begin{enumerate}
	\item assuming that the Riemann Hypothesis is true, prove that $L(s)$ is continuous,
	\item assuming that the Riemann Hypothesis is not true, prove that $L(s)$ cannot be continuous. 
\end{enumerate}
Starting with the first logical step, the hypothesis $\eta(s)\neq 0$ implies that for sufficiently large $n$ also $\eta_n(s)\neq 0$, so that the limit of the following ratio is eventually properly defined. Moreover, because by definition both $\eta_n(s)$ and $\eta_n(1-s)$ converge point wise to $\eta(s)\neq0$ and $\eta(1-s)\neq0$:       
\begin{eqnarray}
	L(s) = \lim_{n\to\infty} X^{\pm}_n(s) =
	\lim_{n\to\infty} 
	\frac{\eta_n(1-s)}{\eta_n(s)} = \frac{\lim_{n\to\infty}\eta_n\left(1-s\right)}{\lim_{n\to\infty}\eta_n\left(s\right)} =
	\frac{\eta(1-s)}{\eta(s)} = X^{\pm}(s) \;\;, 
	\label{eq_19}
\end{eqnarray} 
which is a continuous function on all $D$\,.

\indent
In order to then prove also the second logical step, let us suppose that there exists $s_0\in D$ such that $\eta(s_0)=0$, and consequently also $\eta(1-s_0)=0$, hence resulting in an undefined ratio. By observing that $\eta(s_0)=0$ implies $\eta_{n-1}(s_0)=-R_n(s_0)$, and similarly so for $\eta_{n-1}(1-s_0)=-R_n(1-s_0)$, both $\neq 0$, switching to their polar forms, and recalling \textbf{Theorem 1}, shows that the limit $L(s_0)$ exists nonetheless ($0<\sigma_0<\frac{1}{2}$): 
\begin{eqnarray}
	L(s_0) = \lim_{n\to\infty} X^{\pm}_n(s_0) =
	\lim_{n\to\infty} \frac{|\eta_{n-1}(1-s_0)|}{|\eta_{n-1}(s_0)|}\, e^{i\Delta_n} =
	\lim_{n\to\infty} \frac{|R_n(1-s_0)|}{|R_n(s_0)|} \, e^{i\Delta_n}  = \lim_{n\to\infty}\; \frac{n^{\sigma_0}}{n^{1-\sigma_0}}\, e^{i\Delta_n}   = 0 \label{eq_20}
\end{eqnarray} 
to then coincide with $X^{\pm}(s)$ at any other neighboring point $s \neq s_0$ (zeros of analytic functions are isolated).  It is known that $X^{\pm}(s) \neq 0$ throughout the interior of the critical strip, with $X^{\pm}(s) = 1$ along the critical line, we thus have that $L(s)$ cannot possibly be continuous at $s=s_0$, instead featuring the removable discontinuity $L(s_0)=0$. \qed  \medskip \\ 
Note how \textbf{Theorem 2} cannot provide any additional information about zeros located along the critical line. For such zeros, (\ref{eq_19}) and (\ref{eq_20}) would both yield $L(s)=1$, with no discontinuities. \\
Moreover, \textbf{Theorem 2} indicates $L(s_0)=0$ as the only possible value at discontinuities, if any. Thus, were it possible to prove, from more general properties of the $\eta_n(s)$ sums, that $\lim_{n\to\infty} X^{\pm}_n(s)\neq 0$ everywhere inside the left half of the critical strip, and in so doing proving that $L(s)$ must hence be continuous, then the Riemann Hypothesis would be confirmed. We may finally observe that were the non-trivial zeros simple, then the first derivatives $\eta'(s_0)$ and $\eta'(1-s_0)$ would both be $\neq0$, thus allowing the application of the L'Hôpital's rule in (\ref{eq_19}) to obtain the nested limit:

\begin{eqnarray}
	\lim_{s\to s_0} \,\, \lim_{n\to\infty} \, X^{\pm}_n(s) =
	\lim_{s\to s_0} 
	\frac{\eta(1-s)}{\eta(s)} = \frac{\eta'(1-s_0)}{\eta'(s_0)} = X^{\pm}(s_0) \;\;, 
	\label{eq_21}
\end{eqnarray}    

while (\ref{eq_20}) may also be written as a nested limit:   
\begin{eqnarray}
	\lim_{n\to \infty} \,\, \lim_{s\to s0} \, X^{\pm}_n(s) =
	\lim_{n\to\infty} \frac{|\eta_n(1-s_0)|}{|\eta_n(s_0)|}\, e^{i\Delta_n} = 0 \label{eq_22}
\end{eqnarray} 
Thus, were the Riemann Hypothesis wrong, then the two limiting operations (\ref{eq_21}) and (\ref{eq_22}) would yield differing results, and they cannot therefore be interchanged. However, were it possible to prove, from more general properties of the $\eta_n(s)$ sums, that the sequence $X^{\pm}_n(s)$ converges uniformly on $D$, then said two limiting operations would instead need to be interchangeable, so that assuming the Riemann Hypothesis wrong would lead to a contradiction. 

\section{Final Remarks}
\label{sec:4}         
The key implications of the existence of the above-introduced Hurwitz zero subsequence $\lbrace z_{2k}(s_0)\rbrace$, where $z_{2k}(s_0) = s_0+\delta_{2k}(s_0)$ and $\delta_{2k}(s_0) = \lambda_{2k}(s_0) + i \tau_{2k}(s_0)$, can be effectively visualized by means of the animation shown in Fig. 5 (which starts automatically in Adobe Acrobat Reader). 
For the example $s_0$ corresponding to the sixth known zero, and restricting to indexes $N \leq 2k \leq N + 2m$, where $N=1500$, $2m=270$ (with $2m$ as introduced in (\ref{eq_15}), which we may now conveniently refer to as the \textit{2$\pi$ angular span}), the corresponding 135 individual elements of the sequence $\lbrace z_{2k}(s_0)\rbrace$ were first meticulously computed using a numerical successive-approximations method aimed at achieving $\eta_{2k-1}(z_{2k}) = 0$, or very nearly so. The animation helps illustrate the behavior of this set of indices within the $2\pi$ angular span:
\begin{itemize}
	\item The ring of black dots represents the values $\eta_{2k-1}(s_0) = -R_{2k}(s_0)$ (since, by definition, $\eta(s_0) = 0$), plotted for indices satisfying $N \leq 2k \leq N + 2m$, which correspond to a total angular span of $2\pi$.
	\item The larger gray dot identifies the specific value $\eta_{2k-1}(s_0)$ that is mapped to the origin when $s_0$ is replaced by $z_{2k} = s_0 + \delta_{2k}$.    
	\item The revolving ring of red dots illustrates how all 135 values $\eta_{2k-1}(s_0)$ are consequently displaced to their respective new positions $\eta_{2k-1}(s_0 + \delta_{2k})$. The large red dot located approximately at the center of the ring of red dots represents $\eta(z_{2k}(s_0))$.  
	\item The animation thus shows how the first $2\pi$ turn of the \textit{Deep Spiral 1500} (i.e. $N=1500$) evolves as the $135$ values $z_{2k}(s_0)$ are successively scanned.		
\end{itemize}
\begin{figure}[t]
	\centering
	\animategraphics[loop, autoplay, width=0.40\linewidth]{6}{nr_}{1}{135}
	\caption{ The Hurwitz zeros subsequence $\lbrace z_{2k}(s_0)\rbrace$ brought to life.}
	\label{fig:5}       
\end{figure}
We can observe the following:
\begin{itemize} 	
	\item As the ring of red dots revolves while the 135 values $z_{2k}(s_0)$ are scanned sequentially, its orientation with respect to the stationary ring of black dots (the $\eta_{2k-1}(s_0)$) remains essentially unchanged. In particular, when the black dot highlighted by the larger overlapping gray dot is moved to the origin, all remaining dots are displaced along lines (for example, the blue one) that are very nearly parallel to the line joining the highlighted black dot to the origin (the black line).
	\item This behavior is not a peculiarity of the specific zero chosen for the example, but is instead of very general applicability to all nontrivial zeros. It admits a simple geometric explanation. Referring again to Fig. 4a, the argument of each black segment (i.e., each $u_k(s_0))$ depends only on $k$ and $t_0$, namely $\operatorname{Arg}(u_k(s_0)) = t_0 \log k$. Since $s_0$ is an accumulation point for the sequence ${z_{2k}(s_0)}$, we can always choose $N$ sufficiently large so that the values ${\tau_{2k}(s_0)}$ are guaranteed to be only a tiny fraction of $t_0$. When we restrict attention to the range of $k$ values covering a full interval $t_0\bigl(\log(N+2m) - \log N\bigr)$, corresponding to a $2\pi$ angular span, it follows that the quantities $\tau_{2k}(s_0)\bigl(\log(N+2m) - \log N\bigr)$ can account for at most a tiny fraction of $2\pi$, and in fact tend to zero as $N \to \infty$. Consequently, when passing from $t_0$ to $t_0 + \tau_{2k}(s_0)$, provided $N$ is sufficiently large so that over a $2\pi$ angular span the $\tau_{2k}(s_0)$ values remain only a tiny fraction of $t_0$, the relative orientation of these segments (and hence of each of the 135 gray turns revolving around the black turn in Fig. 5) can only change by at most a very small amount, vanishingly small as $N \to \infty$. Put differently, as we move along the values $z_{2k}(s_0)$, the change in the argument of each $R_{2k}(z_{2k}(s_0))$ is, to a large extent, dictated by the change in $k$ - that is, by the increment from $\log(2k)$ to $\log(2k+2)$ - rather than by the change from $\tau_{2k}$ to $\tau_{2k+2}$.    
\end{itemize}
By the very definition of $\lbrace z_{2k}(s_0)\rbrace$, it is $R_{2k}(z_{2k}(s_0))=\eta(z_{2k}(s_0))$, while (\ref{eq_10}) asymptotically locates the $R_{2k}(z_{2k}(s_0))$ very nearly at the center of the corresponding red ring. In other words, at each sampling point $s=z_{2k}(s_0)$, $R_{2k}(s)$ coincides with the exact value of $\eta(s)$. Hence, for sufficiently large values of $N$, over any range $N \leq 2k \leq N + 2m$ for which the points $\eta_{2k}(s_0)$ sweep out a full angular span of $2\pi$, we obtain that:  
\begin{itemize} 	
	\item The large red dot $\eta(z_{2k}(s_0))$ makes one complete turn around the origin.
	\item Because the black line segment is parallel (or very nearly so) to the line joining the large red dot to the origin, each black segment likewise executes one full revolution around the origin.
	\item At the same time, the blue line segment makes one complete turn around the chosen point $\eta_{2J-1}(s_0)$.
	\item We can therefore conclude that the peculiar geometry governing the convergence of the partial sums of the Dirichlet $\eta$ function forces the winding number of $\eta\bigl(z_{2k}(s_0)\bigr))$ (the large red dot) around $0$ to coincide with the winding number of $\eta_{2J-1}\bigl(z_{2k}(s_0)\bigr)$ (the revolving end of the blue line segment in the animation) around $\eta_{2J-1}(s_0)$, for any fixed $J$. In other words, the two winding numbers are locked together by purely geometric constraints.       
\end{itemize}
Now suppose that the zero $s_0$ is simple, i.e. $\eta'(s_0) \neq 0$. Expanding $\eta$ in a power series around $s_0$ gives $\eta\bigl(z_{2k}(s_0)\bigr) = 0 + \eta'(s_0)\bigl(z_{2k}(s_0) - s_0\bigr) + \cdots$. Since $\eta'(s_0) \neq 0$, the winding number of the curve $\eta\bigl(z_{2k}(s_0)\bigr)$ around $0$ coincides with the winding number of the curve $z_{2k}(s_0)$ around $s_0$. We may therefore represent the latter as a polygonal simple closed curve with $m$ vertices, denoted by $\gamma_N$. The same conclusion about the winding number of $\gamma_N$ around $s_0$ can also be reached by recalling that the winding number of $\eta\bigl(z_{2k}(s_0)\bigr)$ around $0$ coincides with the winding number of $\eta_{2J-1}\bigl(z_{2k}(s_0)\bigr)$ around $\eta_{2J-1}(s_0)$ (locked together and $=1$). Indeed,  $\eta_{2J-1}(z_{2k}(s_0)) = \eta_{2J-1}(s_0) + \eta'_{2J-1}(s_0)\bigl(z_{2k}(s_0) - s_0\bigr) + \cdots$. Hence, also the winding number of the curve $\eta_{2J-1}\bigl(z_{2k}(s_0)\bigr)$ around $\eta_{2J-1}(s_0)$ must coincide with the winding number of $\gamma_N$. On the other hand, arguing by contradiction, suppose instead that $\eta'(s_0) = 0$. Writing $ \eta(s_0) = \eta_{2J-1}(s_0) + R_{2J}(s_0)$, we have $\eta'(s_0) = \eta'_{2J-1}(s_0) + R'_{2J}(s_0) = 0$, which forces $ \eta'_{2J-1}(s_0) = -R'_{2J}(s_0)$. However, by \eqref{eq_16} this quantity is nonzero, so the winding number of $\gamma_N$ would still have to be $1$, yielding a contradiction. \medskip \\ 
So far we have considered only the values of the relevant functions at the sampling points $s=z_{2k}(s_0)$, but not at intermediate values of $s$ along the sides of polygon $\gamma_N$. Since $\eta(s)$ is analytic, and the density of the points $\eta\bigl(z_{2k}(s_0)\bigr)$ increases with $k$ (both radially and angularly; see the discussion following (\ref{eq_15})), we expect that for $s$ lying on the sides of $\gamma_N$ the value $\eta(s)$ remains very close to both sampled values $\eta\bigl(z_{2k}(s_0)\bigr)$ and $\eta\bigl(z_{2k+2}(s_0)\bigr)$. Moreover, because the distance $ \bigl|\eta\bigl(z_{2k}(s_0)\bigr) - \eta\bigl(z_{2k+2}(s_0)\bigr)\bigr|$ is much smaller than $\bigl|\eta\bigl(z_{2k}(s_0)\bigr)\bigr|$, no matter how $\eta(s)$ may wander as $s$ varies along one side of $\gamma_N$, it cannot deviate far enough to encircle $0$ and thus change the winding number. These results are consistent with the conjecture that all nontrivial zeros of $\zeta(s)$ are simple and provide a geometric mechanism supporting this conjecture.

\end{document}